\documentclass[11pt]{amsart}

\usepackage{amssymb, amsmath, amsthm}

 \newtheorem{thm}{Theorem}[section]
 \newtheorem{defn}[thm]{Definition}

 \newcommand{\spec}{\text{spec }}
 \newcommand{\Pspec}{\text{pspec }}
 \newcommand{\Pprim}{\text{pprim }}
 \newcommand{\prim}{\text{prim }}

\begin{document}

\title[natural map]{A natural map from a quantized space onto its semiclassical limit}

\author{Sei-Qwon Oh}

\address{Department of Mathematics, Chungnam National  University, 99 Daehak-ro,   Yuseong-gu, Daejeon 34134, Korea}
 \email{sqoh@cnu.ac.kr}






\begin{abstract}
Let $\widehat{\Gamma}$ be the natural map given in \cite[\S1]{Oh12}. Here, we construct a deformation $B_q$ of a Poisson algebra $B_1$  and a prime ideal $P$ of $B_q$ such that $\widehat{\Gamma}(P)$ is not a Poisson prime ideal of $B_1$.
\end{abstract}

\maketitle



 Recall the star product in \cite[1.1]{Kon}. Let $\hbar$ be a formal variable. Let ${\bf k}$ be a field of characteristic zero and let $A$ be a Poisson algebra over ${\bf k}$.
 The star product $\star$ on $A[[\hbar]]$ is given as follows:
for any  $f=\sum_{n\geq0}f_n\hbar^n, g=\sum_{n\geq0}g_n\hbar^n\in A[[\hbar]]$,
$$f\star g:=\sum_{k,l\geq0}f_kg_l\hbar^{k+l}+\sum_{k,l\geq0,m\geq1} B_m(f_k,g_l)\hbar^{k+l+m}
$$
such that
$$\frac{f\star g-g\star f}{\hbar}\left |\right._{\hbar=0}=\{f_0,g_0\},$$
where $B_i:A\times A\longrightarrow A$ are differential operators with respect to each argument. The ${\bf k}[[\hbar]]$-algebra $A[[\hbar]]$ with a star product is called a {\it quantization} of $A$.

Let $B$ be a ${\bf k}[\hbar]$-subalgebra of $A[[\hbar]]$ containing $A$.  Note that $\hbar$ is a nonzero, nonunit and non-zero-divisor such that $B/\hbar B\cong A$ is commutative.
For $\lambda\in{\bf k}$, $(\hbar-\lambda)B$ is an ideal of $B$ and thus $B_\lambda:=B/(\hbar-\lambda)B$ is a ${\bf k}$-algebra. In particular, we can recover the Poisson algebra $B_0\cong A$ and $B_\lambda=0$ if $\hbar-\lambda$ is a unit in $B$.
Suppose that $0\neq\lambda\in{\bf k}$ such that $\hbar-\lambda$ is a nonunit in $B$. Then  $B_\lambda$ is a non-trivial ${\bf k}$-algebra. (Let $\lambda\in{\bf k}$ be nonzero. Then $\hbar-\lambda$ is a unit in $A[[\hbar]]$. Thus, if $\hbar-\lambda$ is a nonunit in $B$, then $B$ is a proper subalgebra of $A[[\hbar]]$.)

 Let $t$ be an indeterminate. In the above paragraphes, replace
${\bf k}[[h]]$ by ${\bf k}[t,t^{-1}]$ and let $B$ be any ${\bf k}[t,t^{-1}]$-algebra.
Assume that $t-1$ is a nonzero, nonunit and non-zero-divisor in $B$ such that $B_1:=B/(t-1)B$ is commutative as $\hbar$ is such an element in $A[[\hbar]]$.
It is easy to check that $B_1$ is a Poisson algebra with Poisson bracket
\begin{equation}\label{PPPB}
\{\overline{a},\overline{b}\}=\overline{(t-1)^{-1}(ab-ba)}
\end{equation}
for all $\overline{a}=a+(t-1)B,\overline{b}=b+(t-1)B\in B_1$.
The Poisson algebra $B_1$ plays a role as the Poisson algebra $A$ in the first paragraph.
 The Poisson algebra $B_1$ is called a {\it semiclassical limit} of $B$.
 For $\lambda\in{\bf k}$, set $B_\lambda:=B/(t-\lambda)B$ and, for any ${\bf k}$-algebra $R$,  denote by GKdim($R$) the Gelfand-Kirillov dimension of $R$.
 (Refer to \cite{KrLe} for the Gelfand-Kirillov dimension.)
 The ${\bf k}$-algebra $B_\lambda$ is  called a {\it deformation} of $B_1$.
Let   ${\Bbb K}$ be a nonempty subset of the set
$$\{1,0\neq\lambda\in{\bf k}\ |\  t-\lambda \text{ is a nonunit  in $B$ and GKdim($B_\lambda$)=GKdim($B_1$)}\}.$$
That is, if $\lambda\in\Bbb K$ then $B_\lambda\neq0$ and GKdim($B_\lambda$) is equal to GKdim($B_1$).
This article concerns  with the natural map from
a family  of deformations of $B_1$, $\{B_\lambda|\lambda\in\Bbb K\}$,   into the Poisson algebra $B_1$  given in \cite[\S1]{Oh12}.
Recall the natural map given in \cite[\S1]{Oh12}.
 Define a natural homomorphism of ${\bf k}$-algebras
$$\gamma: B\longrightarrow \prod_{\lambda\in{\Bbb K}}B_\lambda,\ \ b\mapsto\gamma(b):=(b+(t-\lambda)B)_{\lambda\in{\Bbb K}}.$$
Let   $\gamma_1:B\longrightarrow B_1$ be the natural projection given by $\gamma_1(b)=b+(t-1)B$ for all $b\in B$.
If  $\gamma$ is a monomorphism, then there exists  the composition
$$\Gamma=\gamma_1\circ\gamma^{-1}:\gamma(B)\longrightarrow  B_1.$$
 Note that $\Gamma$ is surjective.
Denote by $q$ a parameter taking values in ${\Bbb K}$. Then the ${\bf k}$-algebra $B_q$ presents a family of  deformations of $B_1$, $\{B_\lambda| \lambda\in{\Bbb K}\},$  and there is a natural map
$$\widehat{}\ : B_q\longrightarrow \prod_{\lambda\in{\Bbb K}}B_\lambda,\ \ \ z(q)\mapsto(z(\lambda))_{\lambda\in{\Bbb K}}$$
such that the image $\widehat{q}$ of $q$ by the map \ $\widehat{}$ \  is equal to $\gamma(t)$. Note that the image $\widehat{B_q}$ of $B_q$ by the map \ $\widehat{}$ \  is contained in $\gamma(B)$. Hence there is a composition
$$\widehat{\Gamma}=\ \widehat{}\ \circ\Gamma:B_q\longrightarrow B_1,$$
which is the natural map given in \cite[\S1]{Oh12}.

  Deformations $B_q$ of $B_1$ share algebraic properties with the Poisson algebra $B_1$ in many cases, for instance, prime ideals and Poisson prime ideals, primitive ideals and Poisson primitive ideals, central elements and Poisson central elements, automorphisms and Poisson automorphisms. (See \cite{BeKo}, \cite{ChOh3}, \cite{Good4}, \cite{JoOh2}, \cite{Oh7}.)
Let $I$ be an ideal of $B_q$.
  It is proved in  \cite[Theorem 1.4]{Oh12} that  the image $\widehat{\Gamma}(I)$ of $I$ by the map \ $\widehat{\Gamma}$ \ is  a Poisson ideal of $B_1$.
Moreover, it is shown in \cite{Oh12} under certain conditions that the natural map $\widehat{\Gamma}$ induces a homeomorphism between the spectrum of $B_q$ and the Poisson spectrum of $B_1$ for the case of the quantized affine algebra and the quantized Weyl algebra.

The aim of this article is to give an example such that the natural map $\widehat{\Gamma}$ does not preserves primeness. That is, we construct  deformations $B_q$ of a Poisson algebra $B_1$  and a prime ideal $P$ of $B_q$ such that $\widehat{\Gamma}(P)$ is not a Poisson prime ideal of $B_1$.

\bigskip

Assume throughout the article that  ${\bf k}$ denotes a field of characteristic zero and that all algebras  have  unity.
 A commutative algebra $A$  is said to be  a {\it Poisson algebra} if there exists a bilinear product $\{-,-\}$ on $A$, called a {\it Poisson bracket}, such that $(A, \{-,-\})$ is a Lie algebra and $\{ab,c\}=a\{b,c\}+\{a,c\}b$ for all $a,b,c\in A$. An ideal $I$ of a Poisson algebra $A$ is said to be a {\it Poisson ideal} if $\{I,A\}\subseteq I$.  A Poisson ideal $P$ is said to be {\it Poisson prime} if, for all Poisson ideals $I$ and $J$, $IJ\subseteq P$ implies $I\subseteq P$ or $J\subseteq P$.  If $A$ is noetherian  then a Poisson prime ideal of $A$ is a prime ideal by \cite[Lemma 1.1(d)]{Good3}.

\bigskip

Let $B$ be the ${\bf k}[t,t^{-1}]$-algebra generated by $e, f,h$ with the relations
$$ef-fe=(t-1)h,\ \ he-eh=2(t-1)e,\ \ hf-fh=-2(t-1)f.$$
Note that $B$ is an iterated skew polynomial ring
$$B={\bf k}[t,t^{-1}][h][e;\alpha_1][f;\alpha_2,\delta_2]$$
for suitable maps $\alpha_1,\alpha_2, \delta_2$.
(Refer to \cite[Chapter 2]{GoWa2} for the skew polynomial ring.)
Hence  $t-1$ is a non-zero-divisor in $B$. Clearly, $t-1$ is a nonzero and  nonunit such that  $B_1:=B/(t-1)B$ is commutative and thus the semiclassical limit $B_1$ of $B$ is the Poisson algebra ${\bf k}[e,f,h]$ with Poisson bracket
\begin{equation}\label{PB}
\{e,f\}=h,\ \ \{h,e\}=2e,\ \ \ \{h,f\}=-2f
\end{equation}
by (\ref{PPPB}).
Set
 $${\Bbb K}:={\bf k}\setminus(\{0,1\}\cup\{\text{roots of unity}\}).$$  Then ${\Bbb K}$ is  an infinite set clearly and $t-\lambda$, $\lambda\in {\Bbb K}$, is a nonunit in $B$.
 Denote by $B_\lambda$ the deformation $B/(t-\lambda)B$ of $B_1$.
Note that, for each $\lambda\in {\Bbb K}$,
 $B_\lambda$ is  the ${\bf k}$-algebra generated by $e,f,h$ with the relations
\begin{equation}\label{RLAMB}
ef-fe=(\lambda-1)h,\ \ he-eh=2(\lambda-1)e,\ \ hf-fh=-2(\lambda-1)f
\end{equation}
and observe that $B_\lambda$ is an iterated skew polynomial ring $$B_\lambda={\bf k}[h][e;\alpha_1'][f;\alpha_2',\delta_2']$$
for suitable maps $\alpha_1',\alpha_2', \delta_2'$. Hence
 the set of standard monomials $$\{e^if^jh^k|i,j,k=0,1,\cdots\}$$
 forms a ${\bf k}[t,t^{-1}]$-basis of $B$ and a ${\bf k}$-basis of $B_\lambda$ for each $\lambda\in{\bf K}$. It is also a ${\bf k}$-basis of $B_1$.
 Hence GKdim($B_\lambda$)$=$GKdim($B_1$)$=3$ for each $\lambda\in\Bbb K$ by \cite[Example 3.6]{KrLe}.
By  \cite[Lemma 1.2]{Oh12},   there exists a monomorphism of ${\bf k}$-algebras
$$\gamma:B\longrightarrow\prod_{\lambda\in{\Bbb K}}B_\lambda$$
and thus there exists the composition
$$\Gamma=\gamma_1\circ\gamma^{-1}:\gamma(B)\longrightarrow  B_1,$$
where $\gamma_1$ is the canonical projection
$$\gamma_1:B\longrightarrow B_1,\ \ b\mapsto b+(t-1)B.$$
Note that $\Gamma$ is an epimorphism of ${\bf k}$-algebras.

Let $q$ be not a root of unity. Then  $B_q$ is  the ${\bf k}$-algebra generated by $e,f,h$ with the relations
\begin{equation}\label{RQ}
ef-fe=(q-1)h,\ \ he-eh=2(q-1)e,\ \ hf-fh=-2(q-1)f,
\end{equation}
which is obtained from (\ref{RLAMB}) by replacing $\lambda$ by $q$. Moreover,
$q$ plays a role as a parameter taking values in ${\Bbb K}$ and thus
$B_q$ presents a family of deformations of $B_1$, $\{B_\lambda|\lambda\in\Bbb K\}$. There exists a map
$$\widehat{}\ : B_q\longrightarrow \prod_{\lambda\in{\Bbb K}}B_\lambda,\ \ \ z(q)\mapsto(z(\lambda))_{\lambda\in{\Bbb K}}.$$
Observe that the image $\widehat{B_q}$ of $B_q$ by \ $\widehat{}$ \ is contained in $\gamma(B)$. Thus there exists the composition
$$\widehat{\Gamma}=\ \widehat{}\ \circ\Gamma:B_q\longrightarrow B_1.$$
Note that $\widehat{\Gamma}$ is a surjective map but not a ${\bf k}$-linear map since $\widehat{\Gamma}(q)=1$.

The universal enveloping algebra $U(\frak{sl}_2({\bf k}))$ of the Lie algebra $\frak{sl}_2({\bf k})$ is the ${\bf k}$-algebra generated by $E,F,H$ subject to the relations
\begin{equation}\label{RU}
EF-FE=H,\ \ HE-EH=2E,\ \ HF-FH=-2F.
\end{equation}
One should compare (\ref{PB}) with (\ref{RU}).
The map from $U(\frak{sl}_2({\bf k}))$ onto $B_q$ defined by
\begin{equation}\label{ISOM}
E\mapsto (q-1)^{-1}e,\ \ \ F\mapsto (q-1)^{-1}f,\ \ \  H\mapsto (q-1)^{-1}h
\end{equation}
is an isomorphism of algebras by (\ref{RQ}) and (\ref{RU}). Hence $U(\frak{sl}_2({\bf k}))\cong B_q$.

Note that the element  $\Omega=4ef +h^2-2(q-1)h$ of  $B_q$ is a central element of $B_q$ which corresponds to a scalar multiple of the Casimir element of $U(\frak{sl}_2({\bf k}))$.
For each integer $n\geq2$, the ideal $$P_n=\langle e^n, \Omega-(q-1)^2(n^2-1)\rangle$$ of $B_q$ is a prime ideal by (\ref{ISOM}) and \cite[Proposition 4.3]{Cat}. But $\widehat{\Gamma}(P_n)$ is the Poisson ideal
$$\langle e^n, 4ef+h^2\rangle$$ of $B_1$  by \cite[Theorem 1.4]{Oh12},  which is not prime since $n\geq2$.

\bigskip
\noindent
{\bf Acknowledgments}
The author is supported by National  Research Foundation of Korea, Grant 2012-007347.



\bibliographystyle{amsplain}



\providecommand{\bysame}{\leavevmode\hbox to3em{\hrulefill}\thinspace}
\providecommand{\MR}{\relax\ifhmode\unskip\space\fi MR }
\providecommand{\MRhref}[2]{%
  \href{http://www.ams.org/mathscinet-getitem?mr=#1}{#2}
}
\providecommand{\href}[2]{#2}

\end{document}